\documentclass[11pt,twoside]{article}
\usepackage{mathrsfs}
\usepackage{graphics}
\usepackage{amssymb,amsmath}
\usepackage{graphicx}
\usepackage{amsmath}
\usepackage{amsthm}
\usepackage{amsfonts}
\usepackage{amssymb}
\usepackage{latexsym}
\usepackage{amscd,amsxtra}
\usepackage{subfigure}
\usepackage{amscd}
\usepackage{slashed}
\usepackage{titlesec}          

\renewcommand\thesection{\arabic{section}}
\titleformat{\section}{\normalfont\large\bfseries}{\thesection.}{0.5em}{}

\makeatletter  \evensidemargin0.4in \oddsidemargin0.4in

\newtheorem{theorem}{Theorem}[section]
\newtheorem{proposition}{Proposition}[section]
\newtheorem{mainproposition}{Main Proposition}[section]
\newtheorem{lemma}{Lemma}[section]

\numberwithin{equation}{section}

\theoremstyle{definition}
\newtheorem{definition}{Definition}[section]

\theoremstyle{remark}
\newtheorem{remark}{Remark}[section]

\def\ba{\begin{array}}
\def\ea{\end{array}}
\def\be{\begin{equation}}
\def\ee{\end{equation}}
\def\bee{\begin{eqnarray}}
\def\beee{\begin{eqnarray*}}
\def\eee{\end{eqnarray}}
\def\eeee{\end{eqnarray*}}
\def\nn{\nonumber}

\title{\bf  Harmonic maps from degenerating Riemann surfaces\footnotetext{\emph{Date}: \today.}
}
\author{Miaomiao Zhu \thanks{Supported by IMPRS ``Mathematics in
the Sciences'' and the Klaus Tschira Foundation}}
\date{}

\begin{document}
\maketitle

\thispagestyle{empty} \setcounter{page}{1}

\begin{abstract}\vskip 3mm\footnotesize

\noindent We study harmonic maps from degenerating Riemann surfaces
with uniformly bounded energy and show the so-called generalized
energy identity. We find conditions that are both necessary and
sufficient for the compactness in $W^{1,2}$ and $C^{0}$ modulo
bubbles of sequences of such maps.

\vskip 4.5mm




\end{abstract}

\vskip 0.5cm

\section{Introduction} \label{section 1}


Consider a sequence of harmonic maps from compact Riemann surfaces
$(\Sigma_{n},h_{n})$ to a compact Riemannian manifold $(N,g)$,
\bee \label{1.1}u_{n} : \Sigma_{n} \rightarrow {N}, \eee
with uniformly bounded energy $E(u_{n},
\Sigma_{n})\leq\Lambda<\infty$.

In this paper, we study the compactness of the sequence \eqref{1.1}.
We shall first review some well-established analytic aspects related
to this problem and then focus on the case that the domains
$\Sigma_{n}$ degenerate. Our results indicate that when the
topological type of the degeneration is fixed, one can associate to
$(u_{n},\Sigma_{n})$ a sequence of quantities that characterize the
asymptotic behaviour of maps in the limit process.

At first, let us consider the case that the domain surface is fixed,
namely, $\Sigma_{n}=\Sigma$. The uniform energy bound
$E(u_{n})\leq\Lambda$ allows us to find a map $u: \Sigma \rightarrow
N$ such that $u_{n}$ subconverges weakly to $u$. However, in
general, strong convergence fails because of energy concentration at
finitely many points on $\Sigma$, which are called blow-up points
\cite{SU1}, \cite{SU2}. Away from these points, the convergence is
strong. At these points, the ``bubbling" phenomenon can occur and
the concentrated energy can be captured by finitely many bubbles,
i.e., non-trivial harmonic maps from 2-spheres \cite{SU1},
\cite{SU2}. During the ``bubbling" process, there are some necks
joining the base $u: \Sigma \rightarrow N$ to the bubbles or one
bubble to the next. Jost \cite{J} proved that in the limit, these
necks contain no energy, which means all concentrated energy is
captured by the bubbles. Parker \cite{Pa} showed that these necks
actually converge to points in the target manifold, which means that
in the limit the base and the bubbles are connected.

Next, we allow the complex structure on the domain surface to vary.
In this case, consideration of the degeneration of conformal
structures on a Riemann surface will be necessary. Topologically,
the limit surface is obtained by collapsing finitely many simple
closed curves in $\Sigma$. In the end we obtain a surface with nodes
as singularities. There are two types of collapsing curves. The
first is a homotopically trivial one, which corresponds to the
``bubbling" near isolated singularities, for the complex structure
varies in a compact region of the moduli space. The second is a
homotopically nontrivial curve, which corresponds to the
degeneration of the complex structure. By following the ``bubbling"
procedure, we can also find a limit map consisting of a union of
smooth harmonic maps. However, in general, energy may get lost from
some necks, and those necks will fail to converge to points, as in
the explicit example given in \cite{Pa}.

It is worth mentioning that if, in addition, $u_{n}$ are conformal,
i.e., minimal surfaces, then by the technical tools in minimal
surface theory (e.g., the isoperimetric inequality and the
monotonicity property, etc.), we know that in the limit, there is no
energy loss and there are no necks, which gives satisfactory
compactness results. For more details see, for instance, \cite{CS},
\cite{J}, \cite{PW}, or \cite{Pa}. An analogue is Gromov's
compactness theorem for pseudo-holomorphic curves \cite{Gr},
\cite{Wl}, \cite{PW}, \cite{Ye}, \cite{Hu}.

By an asymptotic analysis of harmonic maps from long cylinders,
Chen-Tian \cite{CT} observed that all connecting necks converge
exponentially to geodesics in the target manifold. Moreover, they
\cite{CT} showed that if, in addition, $u_{n}$ is an
energy-minimizing sequence in the same homotopy class, then the
limit geodesics are all of finite length and they contain no energy.
In fact, nontrivial geodesics arise only from the degeneration of
conformal structures, not from ``bubbling".

In order to understand how energy is lost when the domain surfaces
degenerate, we shall give a precise expression of the energy loss
near the nodes and then show the so-called generalized energy
identity for harmonic maps from degenerating surfaces.

Let $u$ be a harmonic map defined on a standard cylinder $P=
[t_{1},t_{2}]\times S^{1}$ with flat metric
$ds^{2}=dt^{2}+d\theta^{2}$. Then the Hopf quadratic differential of
$u$ on $P$ is given by $\Phi(u)=\phi(u)(dt+id\theta)^{2}$, where
\be \label{1.2}\phi(u)= |u_{t}|^2-|u_{\theta}|^{2}-2iu_{t}\cdot
u_{\theta}. \ee
Our main observation is that the following integral:
\be \label{1.3}\int_{\{t\}\times S^{1}}\phi(u) d\theta  \ee
is independent of  $t\in [t_{1},t_{2}]$. Thus \eqref{1.3} defines a
complex number, which we denote by $\alpha=\alpha(u,P)$. We will see
that this quantity can be applied to study the asymptotic behaviour
of the necks appearing near the nodes (or punctures).

Now we consider a sequence of harmonic maps
\bee u_{n}: (\Sigma_{n},h_{n},c_{n}) \rightarrow N  \eee
with uniformly bounded energy $E(u_{n},
\Sigma_{n})\leq\Lambda<\infty$, where $(\Sigma_{n},h_{n},c_{n})$ is
a sequence of closed hyperbolic Riemann surfaces of genus $g>1$ with
hyperbolic metrics $h_{n}$ and compatible complex structures
$c_{n}$. Assume that $(\Sigma_{n},h_{n},c_{n})$ degenerates to a
hyperbolic Riemann surface $(\Sigma,h,c)$ by collapsing $p$ ($1\leq
p\leq 3g-3$) pairwise disjoint simple closed geodesics
$\gamma_{n}^{j},j=1,2,...,p$. For each $j$, the geodesics
$\gamma_{n}^{j}$ degenerate into a pair of punctures
$(\mathcal{E}^{j,1},\mathcal{E}^{j,2})$. Denote the $h_{n}$-length
of $\gamma_{n}^{j}$ by $l^{j}_{n}$, and let $P_{n}^{j}$ be the
standard cylindrical collar about $\gamma_{n}^{j}$.

We associate to the sequence $(u_{n}, \Sigma_{n})$ a sequence of
$p$-tuples $(\alpha_{n}^{1},...,\alpha_{n}^{p})$, where
$\alpha_{n}^{j}:=\alpha(u_{n}, P_{n}^{j}) \in \mathbb{C}$ are the
quantities defined via \eqref{1.3}. Pulling back the hyperbolic
metrics $h_{n}$ and the compatible complex structures $c_{n}$ by
suitable diffeomorphisms $\Sigma \rightarrow \Sigma_{n}\setminus
\cup_{j=1}^{p}{\gamma_{n}^{j}}$ and passing to a subsequence, we can
think of $(h_{n},c_{n})$ as all living on the limit surface $\Sigma$
and converging in $C_{loc}^{\infty}$ to $(h,c)$. Thus, $u_{n}$
becomes a sequence of harmonic maps defined on
$(\Sigma,h_{n},c_{n})$. Then we will show the following generalized
energy identity for harmonic maps from degenerating surfaces:

\begin{theorem} \label{thm1.1} Let $u_{n}:
(\Sigma_{n},h_{n},c_{n}) \rightarrow N$ be a sequence of harmonic
maps with uniformly bounded energy $E(u_{n},
\Sigma_{n})\leq\Lambda<\infty$, where $(\Sigma_{n},h_{n},c_{n})$ is
a sequence of closed hyperbolic Riemann surfaces of genus $g>1$
degenerating to a hyperbolic Riemann surface $(\Sigma,h,c)$ by
collapsing finitely many pairwise disjoint simple closed geodesics
$\{\gamma^{j}_{n},j=1,2,...,p\}$. Then, after selection of a
subsequence, there exist finitely many blow-up points
$\{x_{1},x_{2},...,x_{I}\}$ which are away from the punctures
$\{(\mathcal{E}^{j,1},\mathcal{E}^{j,2}), j=1,2,...p\}$, and
finitely many harmonic maps
\begin{itemize}
\item[] $u:(\overline{\Sigma},\overline{c})\rightarrow N$,
where $(\overline{\Sigma},\overline{c})$ is the normalization of
$(\Sigma,c)$,

\item[]
$\sigma^{i,l}:S^{2}\rightarrow N, l=1,2,...,L_{i}$, near the $i$-th
blow up point $x_{i}$,

\item[]
$ \omega^{j,k}:S^{2}\rightarrow N, k=1,2,...,K_{j}$, near the $j$-th
pair of punctures $(\mathcal{E}^{j,1},\mathcal{E}^{j,2})$,
\end{itemize}
\noindent such that $u_{n}$ converges to $u$ in $C^{\infty}_{loc}$
on $\Sigma \setminus \{x_{i},i=1,2,...,I\}$ and the following holds
\bee \lim \limits_{n \rightarrow \infty}E(u_{n})&=&
           E(u)+\sum_{i=1}^{I}\sum_{l=1}^{L_{i}}E(\sigma^{i,l}) +
           \sum_{j=1}^{p}\sum_{k=1}^{K_{j}}E(\omega^{j,k})+\sum_{j=1}^{p}\lim
           \limits_{n \rightarrow \infty}|{\rm Re} \alpha_{n}^{j}|\cdot
           \frac{\pi^{2}}{l_{n}^{j}}. \nn \\
           \eee
\end{theorem}

Moreover, for each $j$, there are at most finitely many necks
connecting the base $u$ and the bubbles $\omega^{j,k}$. The sum of
the \emph{average lengths} (see Sect. 3) of those necks is
asymptotically equal to
\bee \sqrt{|{\rm Re} \alpha_{n}^{j}|}\cdot
\frac{\pi^{2}}{l_{n}^{j}}.\eee
In fact, we have

\begin{theorem} \label{thm1.2} Assumptions and notations as in Theorem \ref{thm1.1}.
Then
\begin{itemize}
\item[(1)] $(u_{n},\Sigma_{n})$ subconverge in $W^{1,2}$ modulo bubbles,
i.e., in the limit, the necks contain no energy if and only if
\bee \liminf \limits_{n \rightarrow \infty}|{\rm Re} \alpha
_{n}^{j}|\cdot \frac{\pi^{2}}{l_{n}^{j}}=0, \quad j=1,2,...,p. \eee

\item[(2)] $(u_{n},\Sigma_{n})$ subconverge in $C^{0}$ modulo bubbles,
i.e., in the limit, the images of the necks become points if and
only if
\bee \liminf \limits_{n \rightarrow \infty}\sqrt{|{\rm Re} \alpha
_{n}^{j}|}\cdot \frac{\pi^{2}}{l_{n}^{j}}=0, \quad j=1,2,...,p. \eee
\end{itemize}
\end{theorem}

It is clear from the above theorem that the limits $\liminf
\limits_{n \rightarrow \infty}|{\rm Re} \alpha _{n}^{j}|\cdot
\frac{\pi^{2}}{l_{n}^{j}}, j=1,2,...,p$ are the obstructions for
$(u_{n},\Sigma_{n})$ to subconverge in $W^{1,2}$ modulo bubbles, and
the limits $\liminf \limits_{n \rightarrow \infty}\sqrt{|{\rm Re}
\alpha _{n}^{j}|}\cdot \frac{\pi^{2}}{l_{n}^{j}}, j=1,2,...,p$ are
the obstructions for $(u_{n},\Sigma_{n})$ to subconverge in $C^{0}$
modulo bubbles. For each $j$, the asymptotic behaviour of the necks
appearing near the $j$-th node is characterized by
$\{(\alpha_{n}^{j},l_{n}^{j})\}_{n=1}^{\infty}$, namely
\bee E^{j} \approx |{\rm Re} \alpha_{n}^{j}|\cdot
\frac{\pi^{2}}{l_{n}^{j}}, \qquad L^{j} \approx \sqrt{|{\rm Re}
\alpha _{n}^{j}|}\cdot \frac{\pi^{2}}{l_{n}^{j}}, \eee
where $E^{j}$ is the sum of the energies of the necks and $L^{j}$ is
the sum of the average lengths of the necks. Note that the
quantities $\{(\alpha_{n}^{j},l_{n}^{j})\}_{n\geq1}, j=1,2,...,p$
are defined a priori.


For the asymptotics of the imaginary part of $\alpha_{n}^{j}$, we
have the following:

\begin{proposition} \label{pro1.1} Assumptions and notations as in Theorem \ref{thm1.1}. Then
\bee \limsup \limits_{n \rightarrow \infty}|{\rm Im} \alpha
_{n}^{j}|\cdot \frac{\pi^{2}}{l_{n}^{j}}=0, \quad j=1,2,...,p. \eee
\end{proposition}

When the domain surfaces of \eqref{1.1} are degenerating tori, then
the study of the asymptotics of the necks is simpler because of the
fact that any holomorphic quadratic differential on a torus is a
constant. Some modifications to Parker's example \cite{Pa} can
illustrate the asymptotics mentioned, we refer to \cite{Z}.

Wolf \cite{Wo2} studied the asymptotics of families of harmonic maps
between hyperbolic surfaces where the domain degenerates via
pinching finitely many pairwise disjoint simple closed geodesics. In
this case, the energy of the maps goes to infinity. For the
asymptotics of harmonic maps from surfaces to hyperbolic surfaces or
hyperbolic 3-manifolds, where the surfaces degenerate along a
Teichm\"{u}ller ray, see \cite{Wo1}, \cite{Wo3}, \cite{Mi1},
\cite{Mi2}.

There are various energy identities for sequences of different
approximations of harmonic maps from a fixed surface: for a min-max
sequence by Jost \cite{J}; for Struwe's harmonic map flow and
certain Palais-Smale sequences with uniformly $L^{2}$-bounded
tension field, see \cite{St}, \cite{Q1}, \cite{DT}, \cite{LW},
\cite{W}, \cite{QT}, \cite{To1}; for minimizing sequences of
Sacks-Uhlenbeck approximation of harmonic maps by Chen-Tian
\cite{CT}; for the fourth order approximation of harmonic maps, see
\cite{La}. However, the energy identity for general sequences of
Sacks-Uhlenbeck approximations is still open, a natural question
then is whether a certain generalized energy identity holds. Based
on the observations made by Qing \cite{Q2} and Topping \cite{To2},
one expects a complete understanding of the asymptotic behaviour of
the necks appearing near the finite time singularity of the harmonic
map flow. It would be interesting to ask whether one can associate
to the flow suitable quantities that characterize the asymptotics of
the necks mentioned above in a uniform way. If so, then a
classification of the asymptotics in terms of these quantities is
desirable. Recently, Ding et al. \cite{DLL} introduced a flow for
minimal tori and proved the corresponding energy identity.
Considered as certain Palais-Smale sequences from degenerating tori,
its higher genus generalization is expected.

Now, we briefly outline the remaining parts of the paper. In Sect. 2
we recall some preliminary facts about harmonic maps from surfaces.
In Sect. 3 we develop several analytic properties of harmonic maps
from long cylinders. In Sect. 4 we study harmonic maps from
degenerating surfaces and prove Theorems \ref{thm1.1} and
\ref{thm1.2}.

\vskip 0.2cm

\noindent{\bf Acknowledgements} This paper is developed from a
portion of the author's Ph.D. thesis \cite{Z}. He would like to
thank his advisor, Prof. J\"{u}rgen Jost, for encouragement and
inspiration. He would also like to thank Prof. Michael Struwe for
conversations and advice, and Brian Clarke for discussions and help.

\vskip 0.5cm

\section{Preliminaries} \label{section 2}

Let $(\Sigma,h)$ be a Riemann surface with a metric
$h=\lambda^{2}dzd\overline{z}$ in conformal coordinates $z=x+iy$.
Let $(N,g)$ be a compact Riemannian manifold of dimension $d$, and
let its metric in local coordinates be given by $g_{ij}$, with
Christoffel symbols $\Gamma^{i}_{kl}$.

For $u\in{W}^{1,2}(\Sigma,N)$, the energy of $u$ on $\Sigma$ is
\bee \label{2.1} E(u,\Sigma)
=\frac{1}{2}\int_{\Sigma}g_{ij}(u^{i}_{x}u^{j}_{x}+u^{i}_{y}u^{j}_{y})
dxdy.   \eee
A solution of the corresponding Euler-Lagrange equations
\be \label{2.2} \Delta
u^{i}+\Gamma^{i}_{kl}(u^{k}_{x}u^{l}_{x}+u^{k}_{y}u^{l}_{y})=0,
\quad i=1,...,d, \ee
is called a harmonic map. Note that \eqref{2.1} and \eqref{2.2} are
conformally invariant.

If $u$ is in addition conformal, i.e., if the following holds:
\bee
g_{jk}(u_{x}^{j}u_{x}^{k}-u_{y}^{j}u_{y}^{k}-2iu_{x}^{j}u_{y}^{k})=0,
\nn \eee
then $u$ is called a (parametric) minimal surface in $N$.

If we isometrically embed $N$ into some Euclidian space
$\mathbb{R}^{K}$, then \eqref{2.2} can be written as follows:
\be \label{2.3} -\Delta u = A(u)(\nabla u,\nabla u), \ee
where $A(\cdot,\cdot)$ is the second fundamental form of $N$ in
$\mathbb{R}^{K}$. Any $u\in{W}^{1,2}(\Sigma,N)$ that satisfies
\eqref{2.3} weakly is smooth (\cite{H1}, \cite{H2}, or \cite{Ri} for
a new proof).

For $u\in{W}^{1,2}(\Sigma, N)$, the Hopf quadratic differential
associated to $u$ is defined by $\Phi(u)=\phi(u)dz^{2}$, where
\be \phi(u)=|u_{x}|^{2}-|u_{y}|^{2}-2iu_{x}\cdot u_{y}. \nn \ee
\begin{lemma} \label{lem2.1} ${u}$ harmonic $\Rightarrow$ $\phi(u)$
holomorphic. Also, $\phi(u)\equiv0 \Leftrightarrow u$ is conformal.
\end{lemma}
For a proof, see for instance \cite{J}, Lemma 1.2.2.

We list some analytic facts about two-dimensional harmonic maps
proved in \cite{SU1}.

\begin{theorem} \label{thm2.1}
There exists a constant $\epsilon_{0}>0$ that depends only on $N$
such that
\begin{enumerate}
\item[(1)] {\rm(}$\epsilon$-regularity{\rm)} Let $u: D\rightarrow N $ be a smooth
harmonic map satisfying
\bee E(u, D)=\frac{1}{2}\int_{D}|du|^{2} \leq \epsilon_{0}. \nn \eee
Then
\bee \|du\|_{\widetilde{D},1,p} \leq
C(\widetilde{D},p)\|du\|_{D,0,2}, \nn \eee
$\forall \widetilde{D}\subset\subset D$ and $p>1$, where $D$ is some
regular domain in $\mathbb{R}^{2}$, $\widetilde{D}$ is any regular
subdomain of $D$ and $C(\widetilde{D},p)>1$ is a constant depending
only on $\widetilde{D},p$, and the geometry of $N$.

\item[(2)] {\rm(}Singularity removability{\rm)} Let $u$ be a smooth finite-energy
harmonic map from a punctured disk $D\setminus\{0\}$ to $N$. Then
$u$ extends to a smooth harmonic map from $D$ to $N$.

\item[(3)] Any non-trivial harmonic map $u:S^{2}\rightarrow N$ has energy
$E(u)\geq \epsilon_{0}$.
\end{enumerate}
\end{theorem}

Then, we have the following energy identity theorem \cite{SU1},
\cite{J}, \cite{Pa}.

\begin{theorem} \label{thm2.2} Let $\{h_{n}\}$ be a sequence of
Riemannian metrics on $\Sigma$ converging in $C^{\infty}$ to a
Riemannian metric $h$, and let $u_{n}: (\Sigma,h_{n}) \rightarrow
(N,g)$ be a sequence of $h_{n}$-harmonic maps with uniformly bounded
energy $E(u_{n})\leq\Lambda$. Then there are finitely many blow-up
points $\{x_{1},x_{2},...,x_{I}\}\subset \Sigma$, an $h$-harmonic
map $u: (\Sigma,h) \rightarrow (N,g)$ and finitely many nontrivial
harmonic maps $\sigma^{i,l}:S^{2}\rightarrow N, i=1,2,...,I;
l=1,2,...,L_{i}$, such that after selection of a subsequence,
$u_{n}$ converges in $C^{\infty}_{loc}$ to $u$ on $\Sigma \setminus
\{x_{1},x_{2},...,x_{I}\}$, and the following holds
\bee \lim \limits_{n \rightarrow
\infty}E(u_{n})=E(u)+\sum_{i=1}^{I}\sum_{l=1}^{L_{i}}E(\sigma^{i,l}).
\eee
\end{theorem}

During the blow-up process, some necks connecting one bubble to the
next or connecting the base to a bubble appear. Theorem \ref{thm2.2}
shows that in the limit those necks contain no energy. In this case,
we say $u_{n}$ subconverges to $u$ in $W^{1,2}$ modulo bubbles.
Moreover, Parker \cite{Pa} proved that all necks converge to points
in the target manifold, i.e., $u_{n}$ subconverges to $u$ in $C^{0}$
modulo bubbles. Using our terminology, we simply state Parker's
results as follows:

\begin{theorem} (Bubble tree convergence) Notations and assumptions as in Theorem \ref{thm2.2}. Then, after
selection of a subsequence, $u_{n}$ converges to $u$ in $W^{1,2}\cap
C^{0}$ modulo bubbles.
\end{theorem}

For more details on the construction of the bubble trees, see
\cite{PW}, \cite{Pa}.

\vskip 0.5cm
\section{Harmonic maps from cylinders} \label{section 3}

In this section, we study harmonic maps from cylinders and derive
some analytic properties.

Let $P_{T_{1},T_{2}}=[T_{1}, T_{2}]\times S^{1}$ be a standard
cylinder with metric $ds^{2}=dt^{2}+d\theta^{2}$, here
$S^{1}=\mathbb{R}/2\pi \mathbb{Z}$. Since we will only need to
consider long cylinders, w.l.o.g., we always assume that
$T_{2}-T_{1}>2$. Let $u: P_{T_{1},T_{2}}\rightarrow N$ be a $C^{1}$
map. Denote
\be \Theta(t):= \int_{\{t\}\times S^{1}}|u_{\theta}|^{2}. \nn \ee

The following lemma is a modified version of two lemmas proved in
\cite{Pa} and \cite{LW}. For the reader's convenience, we will give
a proof using arguments from \cite{LW}.

\begin{lemma} \label{lem3.1} There exists $\epsilon_{1}>0$, only depending on
$N$, such that if $u: P_{T_{1},T_{2}}\rightarrow N$ is a harmonic
map and
\be \underset{P_{T_{1},T_{2}}}{\sup}|\nabla u|\leq \epsilon_{1}, \nn
\ee
then
\be \label{3.1} \frac{d^{2}}{dt^{2}}\Theta(t)\geq \Theta(t), \quad
\forall t\in[T_{1}, T_{2}]. \ee
Moreover, we have
\bee  \label{3.2} \int_{T_{1}}^{T_{2}}(\Theta(t))^{\nu}dt \leq
2\frac{((\Theta(T_{1}))^{\nu}+(\Theta(T_{2}))^{\nu})}{\nu}, \quad
\forall \nu \in (0,1]. \eee

\end{lemma}

\noindent\emph{Proof}. By a straightforward calculation as in
\cite{LW}, Lemma 2.1, we have
\bee\frac{d^{2}}{dt^{2}}\int_{S^{1}}|u_{\theta}|^{2}&\geq&
(2-C\epsilon_{1}^{2})\int_{S^{1}}|u_{\theta t}|^{2}+
(\frac{3}{2}-C\epsilon_{1}^{2})\int_{S^{1}}|u_{\theta
\theta}|^{2}-\epsilon_{1}^{2}\int_{S^{1}}|u_{\theta}|^{2}.\nn\eee
Here $C$ is a constant depending only on the geometry of $N$. If we
choose $\epsilon_{1}>0$ small enough, then
\bee \frac{d^{2}}{dt^{2}}\int_{S^{1}}|u_{\theta}|^{2}\geq
\frac{5}{4}\int_{S^{1}}|u_{\theta
\theta}|^{2}-\frac{1}{4}\int_{S^{1}}|u_{\theta }|^{2}\geq
\int_{S^{1}}|u_{\theta }|^{2}.   \nn \eee
Here, in the last step, we used the Poincar\'{e} inequality on
$S^{1}$. This proves \eqref{3.1}.

Let $\tau_{i}= \Theta(T_{i}), i=1, 2$. Then we can solve the
following 2nd order ODE:
\bee \label{3.9} \ddot{\rho}-\rho &=& 0, \quad T_{1}\leq t \leq T_{2},  \nn \\
\label{3.10} \rho(T_{1})&=&\tau_{1},   \nn  \\
\label{3.11} \rho(T_{2})&=&\tau_{2}.   \nn  \eee
and obtain a solution $ \rho(t)=\lambda e^{t}+\mu e^{-t}$, where
\be
\lambda=\frac{(e^{T_{2}}\tau_{2}-e^{T_{1}}\tau_{1})}{e^{2T_{2}}-e^{2T_{1}}},
\quad
\mu=\frac{e^{T_{1}+2T_{2}}\tau_{2}-e^{2T_{1}+T_{2}}\tau_{1}}{e^{2T_{2}}-e^{2T_{1}}}.
  \nn \ee
Applying the maximum principle, we conclude
\bee 0\leq \Theta(t)\leq \rho(t), \quad \forall t \in [T_{1},T_{2}].
 \nn \eee
Note that $T_{2}>T_{1}, \tau_{1}\geq0, \tau_{2}\geq0$, and $
\nu\in(0,1]$. By direct calculation, we have
\bee \int_{T_{1}}^{T_{2}}(\Theta(t))^{\nu}dt &\leq&
|\lambda|^{\nu}\frac{(e^{\nu T_{2}}-e^{\nu T_{1}})}{\nu}+
|\mu|^{\nu}\frac{(e^{-\nu T_{1}}-e^{-\nu
T_{2}})}{\nu}      \nn  \\
&\leq&
2\frac{|e^{T_{2}}\tau_{2}-e^{T_{1}}\tau_{1}|^{\nu}}{(e^{2T_{2}}-e^{2T_{1}})^{\nu}}\cdot\frac{(e^{\nu
T_{2}}-e^{\nu
T_{1}})}{\nu}   \nn \\
&\leq& 2\frac{(( \tau_{1}))^{\nu}+(
\tau_{2})^{\nu})}{\nu}\nn \\
&=& 2\frac{((\Theta(T_{1}))^{\nu}+(\Theta(T_{2}))^{\nu})}{\nu}. \nn
\eee
This gives \eqref{3.2}. We have thus finished the proof. \hfill
$\square$

\vskip 0.2cm

Combining Lemma \ref{lem3.1} and the ``$\epsilon$-regularity"
portion of Theorem \ref{thm2.1}, we have

\begin{lemma} \label{lem3.2} There exist $\epsilon_{2}>0$ and $C>0$, depending only on $N$, such that if
$u$ is a harmonic map from $P_{T_{1}, T_{2}}$ to $N$ and
\be \omega:= \underset{t\in[T_{1},T_{2}-1]}{\sup}\int_{[t,
t+1]\times S^{1}}|du|^{2} \leq \epsilon_{2}, \nn \ee
then
\bee \label{3.3}\int_{T_{1}}^{T_{2}}\Theta(t)dt \leq C\omega, \quad
\int_{T_{1}}^{T_{2}}\sqrt{\Theta(t)}dt \leq C\omega^{\frac{1}{2}}.
\eee
\end{lemma}

\noindent\emph{Proof}. Let $\epsilon_{0}>0$ be the constant in
Theorem \ref{thm2.1} (``$\epsilon$-regularity"), and let
$\epsilon_{1}>0$ be the constant in Lemma \ref{lem3.1}. Let
$\epsilon_{2}=\min
\{\epsilon_{0},(\frac{\epsilon_{1}}{C_{0}})^{2}\}$, where $C_{0}$ is
a constant to be determined later. If
\bee \omega=\underset{t\in[T_{1},T_{2}-1]}{\sup}\int_{[t, t+1]\times
S^{1}}|du|^{2} \leq \epsilon_{2}, \nn \eee
then changing variables by translating in $u$ and using
$\epsilon$-regularity property with $D=[-1,2]\times S^{1}$ and
$\widetilde{D}=[0,1]\times S^{1}$, we have
\bee \underset{[T_{1}+1,T_{2}-1]\times S^{1}}{\sup}|\nabla u| \leq
C_{1}(\underset{t\in[T_{1},T_{2}-1]}{\sup}\int_{[t, t+1]\times
S^{1}}|du|^{2})^{\frac{1}{2}}=C_{1}\omega^{\frac{1}{2}} \leq
\frac{C_{1}}{C_{0}}\epsilon_{1}. \nn \eee
where $C_{1}>0$ is a constant depending only on $N$, but not on
$T_{1},T_{2}$. We take $C_{0}$ to be the constant $C_{1}$ here. Then
we can apply Lemma \ref{lem3.1} to conclude that $\forall \nu \in
(0,1]$,
\bee \label{3.4} \int_{T_{1}+1}^{T_{2}-1}(\Theta(t))^{\nu}dt \leq
2\frac{((\Theta(T_{1}+1))^{\nu}+(\Theta(T_{2}-1))^{\nu})}{\nu} \leq
C\underset{P_{T_{1}+1,T_{2}-1}}{\sup}|\nabla u|^{2\nu} \leq C \omega
^{\nu}.  \eee
On the other hand, by the Cauchy inequality, it is not hard to
verify that
\bee \label{3.5}
\int_{T_{1}}^{T_{1}+1}(\Theta(t))^{\nu}dt+\int_{T_{2}-1}^{T_{2}}(\Theta(t))^{\nu}dt
\leq 2\omega ^{\nu}, \forall \nu \in (0,1].     \eee
\eqref{3.3} follows from combining \eqref{3.4}, \eqref{3.5} and
taking $\nu=1, \frac{1}{2}$. \hfill $\square$

\vskip 0.2cm

\begin{lemma}\label{lem3.3} Let $u: P_{T_{1},T_{2}}\rightarrow N$ be a harmonic map. Then for
$t\in [T_{1},T_{2}]$,
\be \label{3.6}\int_{\{t\}\times S^{1}}\phi(u)d\theta  \ee
is independent of  $t\in [T_{1},T_{2}]$, where
\be  \phi(u)= |u_{t}|^2-|u_{\theta}|^{2}-2iu_{t}\cdot u_{\theta} \nn
\ee
and $\phi(u)(dt+id\theta)^{2}$ is the Hopf quadratic differential of
$u$ on $P_{T_{1},T_{2}}$.
\end{lemma}

\noindent\emph{Proof}. By Lemma \ref{lem2.1}, we know that if $u$ is
harmonic then $\phi(u)$ is holomorphic. Given $t_{1}$ and $t_{2}$
such that $T_{1}\leq t_{1}\leq t_{2}\leq T_{2}$, consider the
rectangle $R$ bounded by $[t_{1},t_{2}]\times \{0\}, \{t_{2}\}\times
[0,2\pi], [t_{2},t_{1}]\times \{{2\pi}\},$ and $\{t_{1}\}\times
[2\pi,0].$ By Cauchy's integral theorem, we have
\bee  \oint_{\partial R } \phi(u) = \int_{R}
\overline{\partial}\phi(u) = 0,   \nn \eee
i.e.,
\bee \int_{\{t_{1}\}\times S^{1}} \phi(u)d\theta =
\int_{\{t_{2}\}\times S^{1}}\phi(u)d\theta. \nn \eee
Hence \eqref{3.6} is independent of $t\in [T_{1},T_{2}]$. \hfill
$\square$

\begin{definition} \label{def3.1} Let $u: P_{T_{1},T_{2}}\rightarrow N$ be harmonic. Then we
define a complex number
\be \label{3.7} \alpha(u,P_{T_{1},T_{2}}):= \int_{\{t\}\times
S^{1}}\phi(u)d\theta \in \mathbb{C} \ee
that is associated to $u$ along the cylinder $P_{T_{1},T_{2}}$.
\end{definition}

\begin{remark} \label{rem3.2} It follows from Lemma \ref{lem3.3} that $\alpha(u,P_{T_{1},T_{2}})$
is well-defined. Moreover, we have $\alpha(u,
P_{t'_{1},t'_{2}})=\alpha(u, P_{t_{1},t_{2}})$, $\forall
t_{1}<t'_{1}<t'_{2}<t_{2}$.
\end{remark}

\begin{definition} Let $u: P_{T_{1},T_{2}}\rightarrow N$ be a $C^{1}$ map. Then we
call
\bee L(u,P_{T_{1},T_{2}}) &:=&
\int_{T_{1}}^{T_{2}}(\int_{0}^{2\pi}|u_{t}|^{2}d\theta)^{\frac{1}{2}}dt.
\nn \eee
the \emph{average length }of $u$ along the cylinder
$P_{T_{1},T_{2}}$.
\end{definition}

\begin{remark} Let $c: [T_{1},T_{2}]\rightarrow N$ be a  $C^{1}$
curve in $N$. Then $u(t,\theta):=c(t)$ is a $\theta$-independent
$C^{1}$ map from $P_{T_{1},T_{2}}$ to $N$. It is easy to verify that
\bee L(u,P_{T_{1},T_{2}})=\sqrt{2\pi}L(c,[T_{1},T_{2}]), \nn \eee
where $L(c,[T_{1},T_{2}])=\int_{T_{1}}^{T_{2}}|\dot{c}(t)|dt$ is the
usual length of the curve $c$.
\end{remark}

\begin{lemma} \label{lem3.4} Let $u: P_{T_{1},T_{2}}\rightarrow N$ be a harmonic map
 with $\alpha=\alpha(u,P_{T_{1},T_{2}})$. Then we have

\noindent(1)
\be \label{3.8} |E(u,P_{T_{1},T_{2}})- \frac{1}{2}|{\rm Re}
\alpha|\cdot(T_{2}-T_{1})| \leq \int_{T_{1}}^{T_{2}}\Theta(t)dt, \ee
\noindent(2)
\be \label{3.9} |L(u,P_{T_{1},T_{2}}) - \sqrt{|{\rm Re}
\alpha|}\cdot(T_{2}-T_{1})| \leq
\int_{T_{1}}^{T_{2}}\sqrt{\Theta(t)}dt, \ee
\noindent(3)
\be \label{3.10} |{\rm Im}
\alpha|\cdot(T_{2}-T_{1})\leq2\sqrt{2E(u,P_{T_{1},T_{2}})}\cdot
\sqrt{\int_{T_{1}}^{T_{2}}\Theta(t)dt}.\ee
\end{lemma}

\noindent\emph{Proof}. In view of Definition \ref{def3.1}, we have
\bee {\rm Re} \alpha = \int_{0}^{2\pi}|u_{t}|^{2}d\theta -
\int_{0}^{2\pi}|u_{\theta}|^{2}d\theta, \quad {\rm Im} \alpha = -
2\int_{0}^{2\pi}u_{t}\cdot u_{\theta}d\theta.  \nn \eee
Then by applying the following inequalities,
\bee |(a+b) - |a|| \leq b, \hspace{10pt}  |\sqrt{a+b} - \sqrt{|a|}|
\leq \sqrt{b}, \quad \forall a,b, a+b \geq 0, b \geq 0, \nn  \eee
and then integrating with respect to $t$, we get \eqref{3.8} and
\eqref{3.9}. \eqref{3.10} follows from the Cauchy inequality. \hfill
$\square$

\vskip 0.2cm

The next lemma is inspired by \cite{Pa}.

\begin{lemma} \label{lem3.5} Let $u: P_{T_{1},T_{2}}\rightarrow N$ be a $C^{1}$ map. Then
\be \label{3.11} \underset{P_{T_{1},T_{2}}}{\rm osc}u \leq 4\pi
\underset {P_{T_{1},T_{2}}}{\sup}|\nabla u|+ \frac{1}{\sqrt{2\pi}}
L(u,P_{T_{1},T_{2}}). \ee
\end{lemma}

\noindent\emph{Proof}. Let $(t_{1},\theta_{1}),(t_{2},\theta_{2})\in
P_{T_{1},T_{2}} = [T_{1},T_{2}]\times S^{1}$, where $T_{1}\leq
t_{1}<t_{2}\leq T_{2}$. Then by the Mean Value Theorem for
integration, there exists $ \theta _{0}\in [0,2\pi]$ such that
\bee \label{3.12} \int_{t_{1}}^{t_{2}}|u_{t}(t,\theta
_{0})|dt=\frac{1}{2\pi} \int_{0}^{2\pi}\int_{t_{1}}^{t_{2}}|u_{t}|dt
d\theta. \eee
Hence, we have
\bee {\rm dist}(u(t_{1},\theta_{1}),u(t_{2},\theta_{2})) &\leq&
{\rm dist}(u(t_{1},\theta_{1}),u(t_{1},\theta_{0})) + {\rm dist}(u(t_{1},\theta_{0}),u(t_{2},\theta_{0})) \nn  \\
&&{} + {\rm dist}(u(t_{2},\theta_{0}),u(t_{2},\theta_{2}))                    \nn  \\
&=& I+II+III, \nn \eee
It is easy to see that $ I+III
 \leq\int_{0}^{2\pi}|u_{\theta}(t_{1},\theta)|d\theta
+\int_{0}^{2\pi}|u_{\theta}(t_{2},\theta)|d\theta \leq4\pi \cdot
\underset{P_{T_{1},T_{2}}}{\sup}|\nabla u|$. By \eqref{3.12} and the
Cauchy inequality, we conclude
\bee II &\leq& \int_{t_{1}}^{t_{2}}|u_{t}(t,\theta_{0})|dt =
\frac{1}{2\pi}\cdot \int_{0}^{2\pi}\int_{t_{1}}^{t_{2}}|u_{t}|dt
d\theta \leq \frac{1}{\sqrt{2\pi}}\cdot L(u,P_{T_{1},T_{2}}). \nn
\eee
\eqref{3.11} follows immediately. \hfill $\square$

\vskip 0.2cm

Based on the neck analysis in \cite{DT}, we have the following
proposition, which gives a refined ``bubble domain and neck domain"
decomposition for a sequence of harmonic maps from long cylinders
under certain assumptions.

\begin{proposition} \label{pro3.1} Let $u_{n}\in C^{\infty}(P_{n},N)$ be a
sequence of harmonic maps with $\alpha_{n}=\alpha(u_{n},P_{n})$,
where $P_{n}=[T_{n}^{1},T_{n}^{2}]\times S^{1}$. Assume that:
\begin{itemize}
\item[(1)] ``Long cylinder property"
\be \label{3.13} 1 \ll T_{n}^{1} \ll T_{n}^{2}, \quad i.e., \lim
_{n\rightarrow \infty} \frac{1}{T_{n}^{1}}= 0, \lim _{n\rightarrow
\infty} \frac{T_{n}^{1}}{T_{n}^{2}} = 0, \ee
\item[(2)] ``Uniform energy bound"
\be \label{3.14} E(u_{n},P_{n}) \leq \Lambda < \infty, \ee
\item[(3)] ``Asymptotic boundary conditions"
\bee \label{3.15} \lim_{{n\rightarrow \infty}}
\omega(u_{n},P_{T_{n}^{1},T_{n}^{1}+R}) &=& \lim_{{n\rightarrow
\infty}} \omega(u_{n},P_{T_{n}^{2}-R,T_{n}^{2}})
= 0, \quad \forall R \geq 1, \\
\lim_{{n\rightarrow \infty}}\underset{P_{T_{n}^{1},T_{n}^{1}+1}}{\rm
osc}\ u_{n} &=& \lim_{{n\rightarrow
\infty}}\underset{P_{T_{n}^{2}-1,T_{n}^{2}}}{\rm osc}\ u_{n} = 0,
\nn \eee
where
\be \omega(u,P_{T_{1},T_{2}})
:=\underset{t\in[T_{1},T_{2}-1]}{\sup}\int_{[t,t+1]\times
S^{1}}|du|^{2}. \nn  \ee
\end{itemize}

\noindent Then, after selection of a subsequence, which we still
denote by $(u_{n},P_{n})$, either
\begin{itemize}
\item[(I)]
\be \label{3.16} \lim_{{n\rightarrow \infty}} \omega(u_{n},P_{n}) =
0, \ee
or
\item[(II)]
$\exists K > 0$ independent of $n$ and $2K$ sequences
$\{a_{n}^{1}\}, \{b_{n}^{1}\}, \{a_{n}^{2}\}, \{b_{n}^{2}\},...,
\{a_{n}^{K}\}, \{b_{n}^{K}\}$  such that
\be \label{3.17} T_{n}^{1} \leq a_{n}^{1} \ll b_{n}^{1} \leq
a_{n}^{2} \ll
 b_{n}^{2} \leq ... \leq a_{n}^{K} \ll b_{n}^{K} \leq T_{n}^{2} \ee
and
\be  \label{3.18} (b_{n}^{i}-a_{n}^{i}) \ll T_{n}^{2}, \quad
i=1,2,...,K. \ee
Denote
\be J_{n}^{j}:=[a_{n}^{j},b_{n}^{j}]\times S^{1}, \quad j=1,2,...,K,
\nn \ee
\be I_{n}^{0}:= [T_{n}^{1}, a_{n}^{1}]\times S^{1}, I_{n}^{K}:= [
b_{n}^{K}, T_{n}^{2}]\times S^{1}, I_{n}^{i}:= [b_{n}^{i},
a_{n}^{i+1}]\times S^{1}, i=1,2,...,K-1.  \nn \ee
Then
\begin{itemize}

\item[(i)] $\forall i=0,1,...K,$  $\lim \limits_{{n\rightarrow \infty}}
\omega(u_{n},I^{i}_{n}) = 0$. The maps $u_{n}: I^{i}_{n} \rightarrow
N$ are necks corresponding to collapsing homotopically nontrivial
curves.

\item[(ii)] $\forall j=1,2,...,K,$ there are finitely many
harmonic maps $\omega^{j,l}: S^{2}\rightarrow N, l=1,2,...,L_{j}$,
such that:
\bee \lim \limits_{n \rightarrow
\infty}E(u_{n},J_{n}^{j})=\sum_{l=1}^{L_{j}}E(\omega^{j,l}). \nn
\eee
\end{itemize}
\end{itemize}
\end{proposition}

\noindent\emph{Proof}. If $ \liminf \limits_{{n\rightarrow
\infty}}\omega(u_{n},P_{n})= 0$, then, after selection of a
subsequence, we get \eqref{3.16}. Otherwise, w.l.o.g., we can assume
that
\be \lim\limits_{{n\rightarrow \infty}}\omega(u_{n},P_{n}) = \lim
\limits_{{n\rightarrow \infty}}
\underset{t\in[T_{n}^{1},T_{n}^{2}-1]}{\sup}\int_{[t,t+1]\times
S^{1}}|d u_{n}|^{2}dtd\theta > 0. \nn \ee
Then $\exists \epsilon >0, \{t_{n}\}\in [T_{n}^{1},T_{n}^{2}-1]$,
such that for all $n$ large enough,
\be  \int_{[t_{n},t_{n}+1]\times S^{1}}|d u_{n}|^{2}dtd\theta \geq
\epsilon. \nn  \ee
It follows from the ``asymptotic boundary conditions" \eqref{3.15}
that $t_{n}-T_{n}^{1} \rightarrow \infty$ and $T_{n}^{2}-t_{n}
\rightarrow \infty$. By translation $t\rightarrow t-t_{n}$, we can
think of $u_{n}$ as a harmonic map defined on $[-R_{n},R_{n}]\times
S^{1}$ with $R_{n}\rightarrow \infty$ and
\be \label{3.19} \int_{[0,1]\times S^{1}}|d u_{n}|^{2}dtd\theta \geq
\epsilon, \quad \int_{[-R_{n},R_{n}]\times S^{1}}|d
u_{n}|^{2}dtd\theta \leq \Lambda. \ee
As $n\rightarrow \infty$, $[-R_{n},R_{n}]\times S^{1}$ exhaust
$(-\infty,\infty)\times S^{1}$, which is conformally equivalent to
$S^{2}$ with two punctures. Hence by the conformal invariance of
two-dimensional harmonic maps and the bubble tree convergence
theorem, we can choose a subsequence of $(u_{n},P_{n})$ (still
denoted by $(u_{n},P_{n})$) such that there exist $\{a_{n}^{1}\},
\{b_{n}^{1}\}$ satisfying
\be T_{n}^{1} \leq a_{n}^{1} \ll b_{n}^{1} \leq T_{n}^{2},\quad
(b_{n}^{1}-a_{n}^{1}) \ll T_{n}^{2} \nn \ee
and such that the following sequence of harmonic maps
\bee \widetilde{u}_{n}^{1}(t,\theta):
[-\frac{b_{n}^{1}-a_{n}^{1}}{2},\frac{b_{n}^{1}-a_{n}^{1}}{2}]\times
S^{1} \rightarrow N \nn \eee
converges to a bubble tree $\widetilde{u}_{\infty}^{1}$ (c.f.
\cite{Pa}), where $\widetilde{u}_{n}^{1}(t,\theta):= u_{n}(t+
\frac{a_{n}^{1}+b_{n}^{1}}{2},\theta)$. Moreover, there exist
finitely many harmonic maps $\omega^{1,l} : S^{2}\rightarrow N,
l=1,2,...,L_{1}$, such that
\bee \label{3.20}\lim \limits_{n \rightarrow
\infty}E(u_{n},[a_{n}^{1},b_{n}^{1}]\times S^{1})=\lim \limits_{n
\rightarrow \infty}E(\widetilde{u}_{n}^{1})
=\sum_{l=1}^{L_{1}}E(\omega^{1,l}). \eee
By \eqref{3.19}, \eqref{3.20} and Theorem \ref{thm2.1}, we have $
\limsup\limits_{n \rightarrow \infty}E(u_{n},P_{n}\setminus
([a_{n}^{1},b_{n}^{1}]\times S^{1}))\leq\Lambda-\epsilon_{0}$.
Denote
\be J_{n}^{1}:=[a_{n}^{1},b_{n}^{1}]\times S^{1}, \quad I_{n}^{0}:=
[T_{n}^{1}, a_{n}^{1}]\times S^{1}, \quad I_{n}^{1}:= [ b_{n}^{1},
T_{n}^{2}]\times S^{1}. \nn \ee
Then \eqref{3.20} becomes $ \lim \limits_{n \rightarrow
\infty}E(u_{n},J_{n}^{1})=\sum_{l=1}^{L_{1}}E(\omega^{1,l}).$ After
selection of a subsequence, $u_{n}:I_{n}^{i}\rightarrow N, i=0,1$,
satisfy the conditions \eqref{3.13}, \eqref{3.14} and \eqref{3.15}
with $\Lambda$ replaced by $\Lambda-\epsilon_{0}$. Now if the
following hold:
\be \liminf \limits_{{n\rightarrow \infty}} \omega(u_{n}, I^{i}_{n})
= 0, \quad i=0,1, \nn \ee
then, after passing to a further subsequence, we finish the proof.
Otherwise, we do the same procedure as in the beginning of the proof
and take subsequences if necessary. The whole procedure ends within
finitely many steps because of the uniform energy bound
\eqref{3.14}. The proof can be completed by induction on $K$, the
number of the bubble trees, and reordering
$\{a_{n}^{i},b_{n}^{i}\},i=1,2,...,K.$ \hfill $\square$

Now we study the limit of the energy and average lengths of the
necks
\bee u_{n}: I_{n}^{i}\rightarrow N, \quad i=0,1,...,K. \nn \eee
Recall that these necks satisfy $\lim \limits_{{n\rightarrow
\infty}} \omega(u_{n},I^{i}_{n}) = 0$, hence we can apply Lemma
\ref{lem3.2} and Lemma \ref{lem3.4} to estimate $E(u_{n},I_{n}^{i})$
and $L(u_{n},I_{n}^{i})$.

\begin{mainproposition} \label{mainpro3.1} Assumptions and
notations as in Proposition \ref{pro3.1}, w.l.o.g., we assume that
both $\lim \limits_{n \rightarrow \infty}|{\rm Re} \alpha
_{n}|\cdot|P_{n}|$ and $\lim \limits_{n \rightarrow
\infty}\sqrt{|{\rm Re} \alpha _{n}|}\cdot|P_{n}|$ exist in
$[0,+\infty]$. Then

\noindent(1)
\be \label{3.21} \lim \limits_{n \rightarrow
\infty}\sum_{i=0}^{K}E(u_{n},I_{n}^{i})= \frac{1}{2}\lim \limits_{n
\rightarrow \infty}|{\rm Re} \alpha _{n}|\cdot|P_{n}|, \ee

\noindent(2)
\be
 \lim  \label{3.22}  \limits_{n \rightarrow
\infty}\sum_{i=0}^{K}L(u_{n},I_{n}^{i})= \lim \limits_{n \rightarrow
\infty}\sqrt{|{\rm Re} \alpha _{n}|}\cdot|P_{n}|.  \ee
\end{mainproposition}

\noindent\emph{Proof}. We write
\bee \label{3.23} \sum_{i=0}^{K}E(u_{n},I_{n}^{i}) &=&
\sum_{i=0}^{K}\frac{1}{2}|{\rm Re} \alpha _{n}|\cdot |I^{i}_{n}| +
\sum_{i=0}^{K}(E(u_{n},I_{n}^{i}) - \frac{1}{2}|{\rm Re} \alpha
_{n}|\cdot |I^{i}_{n}|) \nn \\
&=& I+II,  \eee
where
\bee \label{3.24} I &:=& \sum_{i=0}^{K}\frac{1}{2}|{\rm Re} \alpha _{n}|\cdot |I^{i}_{n}|   \nn  \\
&=& \frac{1}{2}|{\rm Re} \alpha _{n}|\cdot[(T_{n}^{2}-T_{n}^{1})-
\sum_{i=1}^{K}(b_{n}^{i}-a_{n}^{i})]                           \nn  \\
&=& \frac{1}{2}|{\rm Re} \alpha _{n}|\cdot
(T_{n}^{2}-T_{n}^{1})\cdot(\frac{T_{n}^{2}}{T_{n}^{2}-T_{n}^{1}})\cdot
[(1-\frac{T_{n}^{1}}{T_{n}^{2}})-
\sum_{i=1}^{K}\frac{(b_{n}^{i}-a_{n}^{i})}{T_{n}^{2}}] \eee
and
\bee II := \sum_{i=0}^{K}(E(u_{n},I_{n}^{i}) - \frac{1}{2}|{\rm Re}
\alpha _{n}|\cdot |I^{i}_{n}|). \nn \eee
Denote $\Theta_{n}(t)= \int_{\{t\}\times
S^{1}}|(u_{n})_{\theta}|^{2}$. By Lemma \ref{lem3.2}, Lemma
\ref{lem3.4} and Proposition \ref{pro3.1},
\bee \label{3.25} |II| &\leq& \sum_{i=0}^{K}|E(u_{n},I_{n}^{i}) -
\frac{1}{2}|{\rm Re}
\alpha _{n}|\cdot |I^{i}_{n}||       \nn  \\
&\leq& \int_{T_{n}^{1}}^{a_{n}^{1}}\Theta_{n}(t)dt +
\sum_{i=1}^{K-1}\int_{b_{n}^{i}}^{a_{n}^{i+1}}\Theta_{n}(t)dt +
\int_{b_{n}^{K}}^{T_{n}^{2}}\Theta_{n}(t)dt  \nn  \\
&\leq& C(\Lambda)\sum_{i=0}^{K}\omega(u_{n}, I^{i}_{n})\rightarrow
0, n\rightarrow \infty, \eee
We write
\bee \label{3.26} \sum_{i=0}^{K}L(u_{n},I_{n}^{i})&=&
\sum_{i=0}^{K}\sqrt{|{\rm Re} \alpha _{n}|}\cdot |I^{i}_{n}| +
\sum_{i=0}^{K}(L(u_{n},I_{n}^{i}) - \sqrt{|{\rm Re} \alpha
_{n}|}\cdot |I^{i}_{n}|) \nn\\
&=& III+IV,  \eee
where
\bee \label{3.27} III &:=& \sum_{i=0}^{K}\sqrt{|{\rm Re} \alpha _{n}|}\cdot |I^{i}_{n}|     \nn  \\
&=& \sqrt{|{\rm Re} \alpha _{n}|}\cdot[(T_{n}^{2}-T_{n}^{1})-
\sum_{i=1}^{K}(b_{n}^{i}-a_{n}^{i})]                                  \nn  \\
&=& \sqrt{|{\rm Re} \alpha _{n}|} \cdot
(T_{n}^{2}-T_{n}^{1})\cdot(\frac{T_{n}^{2}}{T_{n}^{2}-T_{n}^{1}})\cdot
[(1-\frac{T_{n}^{1}}{T_{n}^{2}})-
\sum_{i=1}^{K}\frac{(b_{n}^{i}-a_{n}^{i})}{T_{n}^{2}}] \eee
and
\bee IV := \sum_{i=0}^{K}(L(u_{n},I_{n}^{i})-\sqrt{|{\rm Re} \alpha
_{n}|}\cdot |I^{i}_{n}|). \nn \eee
Applying Lemma \ref{lem3.2}, Lemma \ref{lem3.4} and Proposition
\ref{pro3.1}, we get
\bee \label{3.28} |IV| &=&
\sum_{i=0}^{K}|L(u_{n},I_{n}^{i})-\sqrt{|{\rm Re} \alpha
_{n}|}\cdot |I^{i}_{n}||                                                    \nn  \\
&\leq& \int_{T_{n}^{1}}^{a_{n}^{1}}\sqrt{\Theta_{n}(t)}dt +
\sum_{i=1}^{K-1}\int_{b_{n}^{i}}^{a_{n}^{i+1}}\sqrt{\Theta_{n}(t)}dt
+\int_{b_{n}^{K}}^{T_{n}^{2}}\sqrt{\Theta_{n}(t)}dt                             \nn  \\
&\leq& C(\Lambda)\sum_{i=0}^{K}\sqrt{\omega(u_{n}, I^{i}_{n})}
\rightarrow 0, n\rightarrow \infty. \eee
Recall the properties \eqref{3.13}, \eqref{3.17} and \eqref{3.18} in
Proposition \ref{pro3.1}, namely
\be 1 \ll T_{n}^{1} \ll T_{n}^{2},  \quad 1 \ll(b_{n}^{i}-a_{n}^{i})
\ll T_{n}^{2}, \hspace{20pt} i=1,2,...,K.  \nn \ee
Then, combining \eqref{3.23}, \eqref{3.24} and \eqref{3.25}, we
conclude
\bee \lim \limits_{n \rightarrow
\infty}\sum_{i=0}^{K}E(u_{n},I_{n}^{i}) = \lim\limits_{n\rightarrow
\infty}(I+II) =\frac{1}{2}\lim \limits_{n \rightarrow \infty}|{\rm
Re} \alpha _{n}|\cdot|P_{n}|. \nn \eee
Similarly, combining \eqref{3.26}, \eqref{3.27} and \eqref{3.28}
gives
\bee \lim \limits_{n \rightarrow
\infty}\sum_{i=0}^{K}L(u_{n},I_{n}^{i}) = \lim\limits_{n\rightarrow
\infty}(III+IV) =\lim \limits_{n \rightarrow \infty}\sqrt{|{\rm Re}
\alpha _{n}|}\cdot|P_{n}|. \nn  \eee
Thus we have proved \eqref{3.21} and \eqref{3.22}.  \hfill $\square$

\begin{remark} It follows from Remark \ref{rem3.2}
that
\bee \alpha(u_{n},I_{n}^{i})=\alpha(u_{n},P_{n}), \quad i=0,1,...,K.
\nn \eee
Thus, we can study the properties of the necks
$u_{n}:I_{n}^{i}\rightarrow N$ in a uniform way, but not separately.
\end{remark}

Applying similar arguments as in the proof of Main Proposition
\ref{mainpro3.1}, we get

\begin{proposition} \label{pro3.2} With the same assumptions and notations
as in
 Proposition \ref{pro3.1}, we have
\be \limsup \limits_{n \rightarrow \infty}|{\rm Re} \alpha
_{n}|\cdot|P_{n}|\leq 2\Lambda, \quad \limsup \limits_{n \rightarrow
\infty}|{\rm Im} \alpha _{n}|\cdot|P_{n}|=0,  \nn \ee
\end{proposition}

\noindent\emph{Proof}. By Lemma \ref{lem3.2}, Lemma \ref{lem3.4},
Proposition \ref{pro3.1} and Main Proposition \ref{mainpro3.1}.
\hfill $\square$

\begin{theorem} \label{thm3.1} Assumptions and notations as in
Proposition \ref{pro3.1}. Then

\begin{itemize}

\item[(1)]  $(u_{n},P_{n})$ subconverge in $W^{1,2}$ modulo bubbles,
i.e., in the limit, the necks contain no energy if and only if
\bee \liminf\limits _{n\rightarrow \infty}|{\rm Re} \alpha
_{n}|\cdot |P_{n}|=0. \eee

\item[(2)] $(u_{n},P_{n})$ subconverge in $C^{0}$ modulo bubbles,
i.e., in the limit, the images of the necks become points if and
only if
\bee \liminf\limits_{n\rightarrow \infty} \sqrt{|{\rm Re} \alpha
_{n}|}\cdot |P_{n}|=0. \eee
\end{itemize}
\end{theorem}

\noindent\emph{Proof}.

\noindent(1)  The result is a direct consequence of the identity
\eqref{3.21} in Main Proposition \ref{mainpro3.1}.

\noindent(2) ``$\Leftarrow$": If $\liminf\limits_{n\rightarrow
\infty} \sqrt{|{\rm Re} \alpha _{n}|}\cdot |P_{n}|=0$, then by Main
Proposition \ref{mainpro3.1} and passing to subsequences if
necessary, we have $ \lim \limits_{{n\rightarrow
\infty}}L(u_{n},I_{n}^{i})=0, i=0,1,...,K.$ On the other hand, by
``$\epsilon$-regularity" and the fact that $\lim
\limits_{{n\rightarrow \infty}} \omega(u_{n},I^{i}_{n}) = 0$,
$i=0,1,...,K, \nn $ we get $ \lim \limits_{{n\rightarrow
\infty}}\underset {I_{n}^{i}}{\sup} |\nabla u_{n}|=0, i=0,1,...,K$.
Here we used the fact that, after passing to subsequences, the local
energy of $u_{n}$ over a neighborhood of the two boundary components
of $I^{i}_{n}$ can be arbitrary small. Finally, applying Lemma
\ref{lem3.5}, we conclude
\bee \sum_{i=0}^{K} \underset{I_{n}^{i}}{\rm osc}\ u_{n}
 \leq \sum_{i=0}^{K} (4\pi \cdot \underset {I_{n}^{i}}{\sup}|\nabla u_{n}|
+ \frac{1}{\sqrt{2\pi}}\cdot L(u_{n},I_{n}^{i}))\rightarrow 0, \quad
n \rightarrow \infty. \nn \eee
Thus, all necks converge to points in the target.

  ``$\Rightarrow$": If $(u_{n},P_{n})$ subconverges in $C^{0}$
modulo bubbles, then by the bubble and neck decomposition in
Proposition \ref{pro3.1} and passing to subsequences if necessary,
we get
\bee \label{3.31} \lim \limits_{{n\rightarrow \infty}}
\underset{I_{n}^{i}}{\rm osc}\ u_{n}=0, \quad i=0,1,...,K. \eee
Hence, we have $u_{n}(I_{n}^{i})\subset B(y_{i},\rho_{i})$ for some
$y_{i}\in N$ with $\rho_{i}<{\rm min}(\frac{\pi}{2\kappa},{\rm
inj}(y_{i}))$, where $\kappa^{2}$ is an upper bound on the sectional
curvature of $N$. Fix $i\in\{0,1,...,K\}$ and write
$I_{n}^{i}=[t_{n}^{1},t_{n}^{2}]\times S^{1}$. Then the universal
cover of $I_{n}^{i}$ is
\bee \widetilde{I}_{n}^{i}=\{(t,\theta)\in
\mathbb{R}^{2},t\in[t_{n}^{1},t_{n}^{2}]\}. \nn \eee
It is clear that $u_{n}:I_{n}^{i}\rightarrow B(y_{i},\rho_{i})$
lifts to a harmonic map
\bee \widetilde{u}_{n}:\widetilde{I}_{n}^{i}\rightarrow
B(y_{i},\rho_{i}). \nn \eee
Applying the interior gradient bound for harmonic maps \cite{JK},
\cite{J}, we get
\bee \label{3.32}|d\widetilde{u}_{n}(x_{0})|\leq c_{0} \underset
{x\in
B(x_{0},R)}{\max}\frac{d(\widetilde{u}_{n}(x),\widetilde{u}_{n}(x_{0}))}{R}
\eee
provided $B(x_{0},R)\subseteq\widetilde{I}_{n}^{i}$, where $c_{0}$
is a constant depending only on $N$. Let
$t_{0}=\frac{t_{n}^{2}+t_{n}^{1}}{2}$ and take $x_{0}\in
\{t_{0}\}\times \mathbb{R}$. Then
$B(x_{0},\frac{t_{n}^{2}-t_{n}^{1}}{2})\subseteq\widetilde{I}_{n}^{i}$,
and by \eqref{3.32}, we have
\bee |d\widetilde{u}_{n}(x_{0})|\leq
\frac{2c_{0}}{t_{n}^{2}-t_{n}^{1}} \cdot \underset{I_{n}^{i}}{\rm
osc}\ \widetilde{u}_{n}. \nn \eee
Hence, for $(t_{0},\theta)\in \{t_{0}\}\times S^{1}$,
\bee |du_{n}(t_{0},\theta)|\leq \frac{2c_{0}}{t_{n}^{2}-t_{n}^{1}}
\cdot\underset{I_{n}^{i}}{\rm osc}\ u_{n}. \nn \eee
It follows from
Lemma \ref{lem3.3} and Definition \ref{def3.1} that
\bee \label{3.33} \sqrt{|{\rm Re}\alpha_{n}|} =
|\int_{\{t_{0}\}\times
S^{1}}|(u_{n})_{t}|^{2}-|(u_{n})_{\theta}|^{2}d\theta|^{\frac{1}{2}}
\leq \frac{2c_{0}}{t_{n}^{2}-t_{n}^{1}}
\cdot\underset{I_{n}^{i}}{\rm osc}\ u_{n}.  \eee
Multiplying by $|I_{n}^{i}|=|t_{n}^{2}-t_{n}^{1}|$ on both sides of
\eqref{3.33} gives
\bee \sqrt{|{\rm Re}\alpha_{n}|}\cdot |I_{n}^{i}| \leq C_{0}\
\underset{I_{n}^{i}}{\rm osc}\ u_{n}, \nn \eee
where $C_{0}$ is a constant depending only on $N$, but not on $i$
and $n$. Summing up the inequalities on $I_{n}^{i}$ and applying
\eqref{3.31}, we get
\bee  \sum_{i=0}^{K}\sqrt{|{\rm Re} \alpha _{n}|}\cdot |I^{i}_{n}|
\leq C_{0}\sum_{i=0}^{K} \underset{I_{n}^{i}}{\rm osc}\
u_{n}\rightarrow 0, \quad n \rightarrow \infty. \nn \eee
We thus conclude from Main Proposition \ref{mainpro3.1} that
\bee  \lim\limits_{n\rightarrow \infty}
\sum_{i=0}^{K}L(u_{n},I_{n}^{i})=\lim\limits_{n\rightarrow \infty}
\sqrt{|{\rm Re} \alpha _{n}|}\cdot |P_{n}|=\lim\limits_{n\rightarrow
\infty}\sum_{i=0}^{K}\sqrt{|{\rm Re} \alpha _{n}|}\cdot |I^{i}_{n}|
=0. \nn \eee
\hfill $\square$

Combining Proposition \ref{pro3.1} and Main Proposition
\ref{mainpro3.1} gives the following:

\begin{theorem} \label{thm3.2} Assumptions and notations as in Main Proposition \ref{mainpro3.1}.
Then there exist finitely many harmonic spheres
$\omega^{j,l}:S^{2}\rightarrow N, j=1,2,...,K; l=1,2,...,L_{j}$,
such that after selection of a subsequence of $(u_{n},P_{n})$, we
have
\be \lim \limits_{n \rightarrow
\infty}E(u_{n},P_{n})=\sum_{j=1}^{K}\sum_{l=1}^{L_{j}}E(\omega^{j,l})+
\frac{1}{2}\lim \limits_{n \rightarrow \infty}|{\rm Re} \alpha
_{n}|\cdot |P_{n}|. \nn \ee
\end{theorem}

\vskip 0.5cm

\section{Harmonic maps from degenerating surfaces}
\label{section 4}

In order to study the compactness of a sequence of harmonic maps
$u_{n}: \Sigma_{n}\rightarrow N$, we need to know how the domain
surface varies. We collect some well-known facts about hyperbolic
Riemann surface theory and refer to \cite{A}, \cite{B} and \cite{Hu}
for more details.

\vskip 0.2cm

\noindent{\bf Hyperbolic Riemann surfaces.} We only consider
surfaces without boundary. A Riemann surface $(\Sigma, c)$ is an
orientable surface with a complex structure $c$. A hyperbolic
surface $(\Sigma, h)$ is an oriented surface with a complete
Riemannian metric $h$ of constant curvature $-1$ having finite area.
The topological type of a surface is determined by its signature
$(g,k)$, where $k$ is the number of punctures and $g$ is the genus
of the surface obtained by adding a point at each puncture. The type
$(g,k)$ is called general if
$$2g+k>2.$$ By the uniformization theorem, every Riemann surface of
general type can be represented as a quotient $\mathbb{H}/\Gamma$,
where $\mathbb{H}$ is the Poincar\'{e} upper half plane and $\Gamma$
is a Fuchsian group. Thus, it inherits a hyperbolic metric, where
the punctures become ends. Conversely, for any hyperbolic surface
$(\Sigma, h)$, the induced complex structure extends uniquely to a
conformal structure on the compact surface obtained by adding a
point at each puncture. In fact, there is a natural one-to-one
correspondence between complex structures and hyperbolic metrics on
surfaces of general type.

Two surfaces $\Sigma, \Sigma'$ of type $(g,k)$ are called equivalent
if there exists a conformal diffeomorphism $\Sigma\rightarrow
\Sigma'$ preserving the punctures (if there are any). The space of
equivalence classes is called the moduli space $\mathcal{M}_{g,k}$
of Riemann surfaces of type $(g,k)$. The moduli space
$\mathcal{M}_{g,k}$ in general has certain singularities and thus
does not admit a $C^{\infty}$-structure. It has a covering space
that is a manifold, namely the corresponding Teichm\"{u}ller space.
To this end, we fix a topological model surface $\Sigma_{0}$ of
genus $g$ with $k$ punctures and then consider marked surfaces
$(\Sigma, f)$, where $\Sigma$ is a Riemann surface of type $(g,k)$,
and $f: \Sigma \rightarrow\Sigma_{0}$ is a homeomorphism preserving
the punctures. Two marked surfaces $(\Sigma, f)$ and $(\Sigma', f')$
are called equivalent if there exists a conformal diffeomorphism
$\Sigma \rightarrow \Sigma'$ homotopic to $f'^{-1}\circ f$. The
space of equivalence classes is called the Teichm\"{u}ller space
$\mathcal{T}_{g,k}$ of Riemann surfaces of type $(g,k)$.

Now we consider closed Riemann surfaces of genus $g>1$. Any such
surface is of general type and it acquires a complete hyperbolic
metric. Thus, we are working on the compactness of a sequence of
harmonic maps whose domain surface $\Sigma$ varies in
$\mathcal{M}_{g}$. Ideally, we hope the domain varies in a compact
region. Unfortunately, the moduli space $\mathcal{M}_{g}$ is
non-compact because the conformal structure on $\Sigma$ can
degenerate. The following lemma \cite{M} shows that the only process
by which the conformal structure on $\Sigma$ can degenerate is the
shrinking of simple closed geodesics on $\Sigma$. We represent
$\Sigma$ as a quotient $\mathbb{H}/\Gamma$.

\begin{lemma} Let $\{\Gamma_{n}\}$ be a sequence of Fuchsian groups
which are isomorphic as abstract groups and with non-singular
compact quotients $\mathbb{H}/\Gamma_{n}$. Suppose the lengths of
simple closed geodesics on $\mathbb{H}/\Gamma_{n}$ are uniformly
bounded from below by a positive constant. Then a subsequence of
$\{\Gamma_{n}\}$ converges to some Fuchsian group $\Gamma$ which is
isomorphic to all $\Gamma_{n}$. The convergence can be interpreted
as the convergence of suitably normalized fundamental regions.
\end{lemma}

The natural way to compactify $\mathcal{M}_{g}$, then, is to allow
the lengths of the geodesics to become zero and thus admit surfaces
with nodes as singularities. Topologically, one cuts the surface at
a collection of finitely many homotopically independent pairwise
disjoint simple closed curves and pinches the cut curves to points.
This yields the Deligne-Mumford compactification
$\overline{\mathcal{M}}_{g}$, whose boundary
$\overline{\mathcal{M}}_{g}\setminus \mathcal{M}_{g}$ consists of
surfaces with nodes \cite{DM}. On $\mathcal{T}_{g}$, one can use
Fenchel-Nielsen coordinates to describe this process and obtain the
corresponding partial compactification $\overline{\mathcal{T}}_{g}$
(c.f. \cite{A}, \cite{B}).

Here, following \cite{Hu}, we describe this process in terms of
hyperbolic surface theory. Let $\Sigma_{0}$ be a topological model
surface and $\mathscr{E}^{J}=\{\gamma^{j},j\in J\}$ a possibly empty
collection of finitely many pairwise disjoint simple closed curves
on $\Sigma_{0}$. Let $\widetilde{\Sigma}$ be the surface obtained
from $\Sigma_{0}$ by pinching all curves $\gamma^{j}$ to points
$\mathcal{E}^{j}$. We remove all $\mathcal{E}^{j}$ from
$\widetilde{\Sigma}$ and place a complete hyperbolic metric $h$ on
the resulting surface $\Sigma=\widetilde{\Sigma}\setminus \cup_{j\in
J}\mathcal{E}^{j}$. For $j\in J$, we denote by
$(\mathcal{E}^{j,1},\mathcal{E}^{j,2})$ a pair of punctures on
$(\Sigma,h)$ corresponding to $\mathcal{E}^{j}$. Denote by
$\overline{\Sigma}$ the surface obtained by adding a point at each
puncture of $\Sigma$. Then the complex structure $c$ on $\Sigma$
that is compatible with the hyperbolic structure $h$ extends to a
complex structure $\overline{c}$ on $\overline{\Sigma}$.
$(\widetilde{\Sigma},h,\overline{c})$ is called a nodal surface.
$(\overline{\Sigma},\overline{c})$ is called the normalization of
$(\widetilde{\Sigma},h,\overline{c})$ or $(\Sigma,h,c)$.
$\overline{\Sigma}$ is a surface of lower topological type.

Let $(\Sigma_{n},h_{n},c_{n})$ be a sequence of closed hyperbolic
Riemann surfaces of genus $g>1$. We say that
$(\Sigma_{n},h_{n},c_{n})$ converges to a nodal surface
$(\widetilde{\Sigma},h,\overline{c})$ or a hyperbolic surface
$(\Sigma,h,c)$, if there exist possibly empty collections
$\mathscr{E}^{J}_{n}=\{\gamma^{j}_{n},j\in J\}$ of finitely many
pairwise disjoint simple closed geodesics on each
$(\Sigma_{n},h_{n},c_{n})$ and continuous maps $\tau_{n}: \Sigma_{n}
\rightarrow \widetilde{\Sigma}$ with
$\tau_{n}(\gamma_{n}^{j})=\mathcal{E}^{j}$ for $j\in J$ and each
$n$, such that:
\begin{itemize}
\item[(1)] The lengths $\ell(\gamma_{n}^{j})=l_{n}^{j} \rightarrow 0$ for all $j\in J$.

\item[(2)] $\tau_{n}: \Sigma_{n}\setminus \cup_{j\in J}\gamma_{n}^{j}\rightarrow \Sigma$ is a diffeomorphism for each $n$.

\item[(3)] $(\tau_{n})_{*}h_{n} \rightarrow h $ in $C_{loc}^{\infty}$ on $\Sigma$.

\item[(4)] $(\tau_{n})_{*}c_{n} \rightarrow c $ in $C_{loc}^{\infty}$ on $\Sigma$.
\end{itemize}

By the thick-thin decomposition of a closed hyperbolic surface of
genus $g>1$, the number of small simple closed geodesics (of lengths
$< {\rm 2arcsinh(1)}$) is bounded by $3g-3$ (cf. \cite{Hu}, Lemma
IV.4.1). Thus, we have $0 \leq|J|\leq 3g-3$. If $|J|>0$, we say
$(\Sigma_{n},h_{n},c_{n})$ degenerates to a nodal surface
$(\widetilde{\Sigma},h,\overline{c})$ or a hyperbolic surface
$(\Sigma,h,c)$. Using our notations, we state the following
proposition and refer to \cite{Hu} for a detailed proof.

\begin{proposition} \label{pro4.1}
Let $(\Sigma_{n},h_{n},c_{n})$ be a sequence of closed hyperbolic
Riemann surfaces of genus $g>1$. Then, after selection of a
subsequence, $(\Sigma_{n},h_{n},c_{n})$ converges to a nodal surface
$(\widetilde{\Sigma},h,\overline{c})$ or a hyperbolic surface
$(\Sigma,h,c)$.
\end{proposition}

Thus, the analysis of the degeneration of hyperbolic surfaces is
reduced to the local behaviour of the pinched geodesics. A
fundamental tool to realize this localization is the following
collar lemma \cite{K}, \cite{Mt}, \cite{Ha}, \cite{R}. We again
represent a closed Riemann surface of genus $g>1$ as a quotient
$\mathbb{H}/\Gamma$.

\begin{lemma} \label{lem4.2} Let $\gamma$ be a simple closed
geodesic of length $\ell(\gamma)=l$ in $\mathbb{H}/\Gamma$. Then
there is a collar of area $\frac{l}{\sinh(\frac{l}{2})}$ around
$\gamma$, i.e., $\mathbb{H}/\Gamma$ contains an isometric copy of
the region
\bee \label{4.1} A=\left\{z=re^{i\phi}\in \mathbb{H}:1\leq r \leq
e^{l}, \arctan(\sinh(\frac{l}{2}))<\phi<\pi -
\arctan(\sinh(\frac{l}{2}))\right\},  \eee
where $\gamma$ corresponds to $\{re^{i\frac{\pi}{2}}\in
\mathbb{H}:1\leq r \leq e^{l}\}$, and the lines $\{r=1\}$,
$\{r=e^{l}\}$ are identified via $z\rightarrow e^{l}z$.
\end{lemma}


This collar neighborhood is a topological cylinder and its geometry
is determined by the length of the core geodesic and is hence
independent of the surface. There are other versions of the collar
in terms of different coordinates, for example, a hyperbolic
cylinder with Fermi coordinates \cite{B}. In view of the results
developed in Sect. 3, we need a standard cylindrical version of the
collar \eqref{4.1}. To this end, we consider the following conformal
transformation:
\bee \label{4.2} re^{i\phi}\rightarrow
(t,\theta)=(\frac{2\pi}{l}\phi,\frac{2\pi}{l}\log r). \eee
Then the collar $A$ in
 Lemma \ref{lem4.2} is isometric to the following cylinder:
\bee \label{4.3}
P=\left\{(t,\theta):\frac{2\pi}{l}\arctan(\sinh(\frac{l}{2})) <t<
\frac{2\pi}{l}(\pi-\arctan(\sinh(\frac{l}{2}))),0\leq \theta \leq
2\pi\right\}  \eee
with metric $ds^{2}=(\frac{l}{2\pi \sin
\frac{lt}{2\pi}})^{2}(dt^{2}+d\theta^{2})$; here $\gamma \subset A$
corresponds to $\{t=\frac{\pi^{2}}{l}\}\subset P$, and the lines
$\{\theta=0\}$,$\{\theta=2\pi\}$ in \eqref{4.3} are identified.

Let ${\rm injrad}(\phi,r)$ be the injectivity radius at the point
$(\phi,r)$ of $A$. Then, by results from hyperbolic trigonometry
(see \cite{Hu}, Example 5.5 or \cite{B}, Chapter 2), one can verify
that
\bee \label{4.4}\sinh({\rm
injrad}(\phi,r))\sin(\phi)=\sinh(\frac{l}{2}),\quad (\phi,r)\in A.
\eee
Hence, applying the isometric transformation \eqref{4.2}, we have
\bee \label{4.5} \sinh({\rm
injrad}(t,\theta))\sin(\frac{lt}{2\pi})=\sinh(\frac{l}{2}),\quad
(t,\theta) \in P, \eee
where ${\rm injrad}(t,\theta)$ is the injectivity radius at the
point $(t,\theta)$ of $P$.

\begin{remark}
\eqref{4.4} and \eqref{4.5} are very useful in that they give
explicit expressions of the injectivity radius in terms of two
different coordinates of the points in the collar.
\end{remark}

\noindent{\bf Generalized energy identity.} Consider a sequence of
harmonic maps
\bee u_{n}: (\Sigma_{n},h_{n},c_{n}) \rightarrow N,\eee
with uniformly bounded energy $E(u_{n},
\Sigma_{n})\leq\Lambda<\infty$, where $(\Sigma_{n},h_{n},c_{n})$ is
a sequence of closed hyperbolic Riemann surfaces of genus $g>1$ with
hyperbolic metrics $h_{n}$ and compatible complex structures
$c_{n}$. We are only interested in the case that degeneration
occurs. Thus, by Proposition \ref{pro4.1}, we can assume that
$(\Sigma_{n},h_{n},c_{n})$ converges to a hyperbolic Riemann surface
$(\Sigma,h,c)$ by collapsing $p$ ($1\leq p\leq 3g-3$) pairwise
disjoint simple closed geodesics $\gamma_{n}^{j},j=1,2,...,p$.
Denote the $h_{n}$-length of $\gamma_{n}^{j}$ by $l^{j}_{n}$. Then
in the degeneration $(n \rightarrow \infty)$, we have $
l_{n}^{j}\rightarrow 0, j=1,2,...,p$. For each $j$, the geodesics
$\gamma_{n}^{j}$ degenerate into a pair of punctures
$(\mathcal{E}^{j,1},\mathcal{E}^{j,2})$.

\vskip 0.2cm

\noindent\emph{Proof of Theorem \ref{thm1.1}}. We first consider the
simpler case that $p=1$ and hence omit the indices $j$. Since
$\lim_{n\rightarrow \infty}l_{n}=0$, w.l.o.g., we assume that
$l_{n}\leq 2 {\rm arcsinh(1)}$ for all $n$. Let $P_{n} $ be the
cylindrical collar about $\gamma_{n}$ given by \eqref{4.3} and let
$\alpha_{n}=\alpha(u_{n},P_{n})$ be the complex number associated to
$u_{n}$ along the collar $P_{n}$ as in Definition \ref{def3.1}.
After passing to a subsequence, we can assume that the limit $$\lim
\limits_{n \rightarrow \infty}|{\rm Re} \alpha_{n}|\cdot
\frac{\pi^{2}}{l_{n} }$$ exists in $[0,\infty]$. We will eventually
see that the limit is finite, since the total energy of $u_{n}$ is
uniformly bounded.

For $0 <\delta < {\rm arcsinh(1)}$, let $\Sigma^{\delta} :=
\{z\in\Sigma, {\rm injrad}(z;h)\geq\delta\}$ be the $\delta$-thick
part of the hyperbolic surface $(\Sigma,h)$. Recall that there are
diffeomorphisms $\tau_{n}: \Sigma_{n}\setminus \gamma_{n}\rightarrow
\Sigma$ such that $((\tau_{n})_{*}h_{n},(\tau_{n})_{*}c_{n})$
converges to $(h,c)$ in $C_{loc}^{\infty}$ on $\Sigma$. Set
\bee  \overline{u}_{n}=(\tau_{n})_{*}u_{n}, \quad
\overline{h}_{n}=(\tau_{n})_{*}h_{n}, \quad
\overline{c}_{n}=(\tau_{n})_{*}c_{n}, \nn \eee
and consider the following sequence of harmonic maps:
\bee \overline{u}_{n}:(\Sigma,
\overline{h}_{n},\overline{c}_{n})\rightarrow N. \nn \eee
Then for each fixed $\delta>0$,
$(\overline{h}_{n},\overline{c}_{n})$ converges to $(h,c)$ in
$C^{\infty}$ on $\Sigma^{\delta}$. Choose a fixed sequence
$\delta_{n} \searrow 0$ such that $\Sigma^{\delta_{n}}$ exhaust
$\Sigma$. Then by Theorem \ref{thm2.2} and a standard diagonal
argument, there exist finitely many blow-up points
$\{x_{1},x_{2},...,x_{I}\}\subset \Sigma$ which are away from the
punctures $(\mathcal{E}^{1},\mathcal{E}^{2})$, finitely many
harmonic maps $\sigma^{i,l}:S^{2}\rightarrow N, l=1,2,...,L_{i}$,
near the $i$-th blow-up point $x_{i}$ and a harmonic map
$u:(\Sigma,h,c) \rightarrow N$ such that, after selection of a
subsequence of $(u_{n},\Sigma_{n})$,
$\overline{u}_{n}:(\Sigma^{\delta_{n}},
\overline{h}_{n},\overline{c}_{n})\rightarrow N$ converges in
$C_{loc}^{\infty}$ to $u$ on $\Sigma$. $u$ extends smoothly to the
normalization $(\overline{\Sigma},\overline{c})$ of $(\Sigma,h,c)$.
Moreover, the following holds:
\bee \label{4.7} \lim \limits_{n \rightarrow
\infty}E(u_{n},\tau_{n}^{-1}(\Sigma^{\delta_{n}}))=\lim \limits_{n
\rightarrow \infty}E(\overline{u}_{n},\Sigma^{\delta_{n}}) =
E(u)+\sum_{i=1}^{I}\sum_{l=1}^{L_{i}}E(\sigma^{i,l}).
 \eee
It should be remarked that the subsequence $(u_{n},\Sigma_{n})$ can
be taken in such a way that $\lim \limits_{{n\rightarrow \infty}}
{\rm osc}_{\partial \Sigma^{\delta_{n}}} \overline{u}_{n} =0$, or
equivalently, $\lim \limits_{{n\rightarrow \infty}} {\rm
osc}_{\partial (\Sigma_{n}\setminus
\tau_{n}^{-1}(\Sigma^{\delta_{n}}))} u_{n} =0$.

To recover the energy concentration at the punctures
$(\mathcal{E}^{1},\mathcal{E}^{2})$, we have to study
$(\overline{u}_{n}, \Sigma \setminus \Sigma^{\delta_{n}})$, or
equivalently, $(u_{n},\Sigma_{n}\setminus
\tau_{n}^{-1}(\Sigma^{\delta_{n}}))$. For each $n$ and $\delta$,
$\Sigma_{n}\setminus \tau_{n}^{-1}(\Sigma^{\delta})$ is not the
$\delta$-thin part of $(\Sigma_{n},h_{n})$. However, we claim that
for fixed small $\delta>0$ and for $n$ sufficiently large,
$\Sigma_{n}\setminus \tau_{n}^{-1}(\Sigma^{\delta})$ is almost the
$\delta$-thin part of $(\Sigma_{n},h_{n})$.

To see this, fix $\delta>0$ small and let $z\in\Sigma$ be a point
satisfying ${\rm injrad}(z;h)=\delta$. Since $(\tau_{n})_{*}h_{n}$
converges to $h$ in $C_{loc}^{\infty}$ on $\Sigma$, then for all
$\delta_{1}, \delta_{2}>0$ such that $\delta_{1}<\delta<\delta_{2}$,
the following holds:
\bee \label{4.8} \delta_{1}<{\rm injrad}(z;(\tau_{n})_{*}h_{n})<
\delta_{2}, \quad \textrm{for all $n$ large enough}. \eee
Recall that for $0 <\delta < {\rm arcsinh(1)}$, the $\delta$-thin
part of a hyperbolic surface is either an annulus or a cusp (c.f.
\cite{Hu}, Proposition IV.4.2). For $n\geq1$ and $\delta \in
[\frac{l_{n} }{2}, {\rm arcsinh(1)}]$, let us see what the
$\delta$-thin part of $(\Sigma_{n},h_{n})$ looks like. Recall that
$P_{n}$ is the cylindrical collar about $\gamma_{n}$. Now, we define
the following $\delta$-subcollars of $P_{n}$
\bee \label{4.9}
P_{n}^{\delta}:=[T_{n}^{1,\delta},T_{n}^{2,\delta}]\times S^{1}
\subseteq P_{n}, \eee
where
\bee \label{4.10} T_{n}^{1,\delta} &=&
\frac{2\pi}{l_{n}}\arcsin(\frac{\sinh(\frac{l_{n}}{2})}{\sinh
\delta}), \quad T_{n}^{2,\delta} = \frac{2\pi^{2}}{l_{n}}-
\frac{2\pi}{l_{n}}\arcsin(\frac{\sinh(\frac{l_{n}}{2})}{\sinh
\delta}). \eee
By  \eqref{4.5}, one can verify that $P_{n}^{\delta}$ is exactly the
$\delta$-thin part of $(\Sigma_{n},h_{n})$, namely
\bee \label{4.11}  P_{n}^{\delta}=\{z\in\Sigma_{n}, {\rm
injrad}(z;h_{n})\leq\delta\}.  \eee
Thus, fix $\delta>0$ small, for all $\delta_{1}, \delta_{2}>0$
satisfying $\frac{l_{n} }{2}<\delta_{1}<\delta<\delta_{2}<{\rm
arcsinh(1)}$, it follows from \eqref{4.8} and \eqref{4.11} that
\bee \label{4.12} P_{n}^{\delta_{1}} \subseteq \Sigma_{n}\setminus
\tau_{n}^{-1}(\Sigma^{\delta}) \subseteq P_{n}^{\delta_{2}}, \quad
\textrm{for all $n$ large enough}.   \eee
If we choose $\delta_{1}, \delta_{2}$ in \eqref{4.12} sufficiently
close to $\delta$, then for $n$ large enough, $\Sigma_{n}\setminus
\tau_{n}^{-1}(\Sigma^{\delta})$ is almost the $\delta$-thin part
$P_{n}^{\delta}$ of $(\Sigma_{n},h_{n})$. Thus we have verified our
claim.

For $\delta>0$ small and for $n$ large enough, denote
\bee    \Omega_{n}^{\delta}:= \{(\Sigma_{n}\setminus
\tau_{n}^{-1}(\Sigma^{\delta}))\setminus P_{n}^{\delta}\} \cup
\{P_{n}^{\delta}\setminus(\Sigma_{n}\setminus
\tau_{n}^{-1}(\Sigma^{\delta}))\}.  \nn \eee
Note that $P_{n}$ are equipped with hyperbolic metrics which are
conformal to $dt^{2}+d\theta^{2}$. By the conformal invariance of
harmonic maps, we can replace the hyperbolic metrics with the metric
$dt^{2}+d\theta^{2}$. Recall that $\lim \limits_{{n\rightarrow
\infty}} {\rm osc}_{\partial (\Sigma_{n}\setminus
\tau_{n}^{-1}(\Sigma^{\delta_{n}}))} u_{n} =0$. By applying
``$\epsilon$-regularity" and taking subsequences, we have
\bee \lim \limits_{{n\rightarrow \infty}} {\rm
osc}_{\Omega_{n}^{\delta}} u_{n} =0, \quad \lim
\limits_{{n\rightarrow \infty}}E(u_{n},\Omega_{n}^{\delta})=0. \nn
\eee
Thus, after passing to further subsequences, we conclude
\bee \label{4.13}\lim \limits_{n \rightarrow
\infty}E(u_{n},\Sigma_{k_{n}}\setminus
\tau_{n}^{-1}(\Sigma^{\delta_{n}}))= \lim \limits_{n \rightarrow
\infty}E(u_{n},P_{n}^{\delta_{n}}). \eee
Now the energy concentration at the punctures is reduced to the
study of $(u_{n},P_{n}^{\delta_{n}})$. In view of \eqref{4.9} and
\eqref{4.10}, by choosing further subsequences of
$(u_{n},\Sigma_{n})$, we have
\bee |P_{n}^{\delta_{n}}| = \frac{2\pi^{2}}{l_{n}}-
\frac{4\pi}{l_{n}}\arcsin(\frac{\sinh(\frac{l_{n}}{2})}{\sinh
\delta_{n}}) = \frac{2\pi^{2}}{l_{n}}(1+o(1)), \quad n \rightarrow
\infty. \nn \eee
To apply Theorem \ref{thm3.2} with domain cylinders
$P_{n}^{\delta_{n}}$, we see that the first two conditions, the
``long cylinder property"\eqref{3.13} and the ``uniform energy
bound"\eqref{3.14}, are satisfied. We need to check the ``asymptotic
boundary conditions"\eqref{3.15}. For any fixed $R\geq0$ and for
fixed small $\delta>0$, denote
\bee
A_{n}^{1}(\delta,R):=[T_{n}^{1,\delta}-1,T_{n}^{1,\delta}+R]\times
S^{1}, \quad
A_{n}^{2}(\delta,R):=[T_{n}^{2,\delta}-R,T_{n}^{2,\delta}+1]\times
S^{1}. \nn \eee
Then by \eqref{4.5}, one can verify that the injectivity radii of
the points in $A_{n}^{i}(\delta,R), i=1,2$, are uniformly bounded
from below by a positive constant as $n \rightarrow \infty$. Hence
the images $\tau_{n}(A_{n}^{i}(\delta,R)),i=1,2$ are uniformly away
from the punctures of $\Sigma$. Moreover, the energies
$\sum_{i}E(u_{n}, A_{n}^{i}(\delta,R))$ can be uniformly controlled
by $E(u,\Sigma\setminus\Sigma^{\delta'})$ (for some
$\delta'>\delta$), which goes to $0$ as $\delta'\rightarrow 0$.
Thus, after passing to subsequences, one can verify the ``asymptotic
boundary conditions". Now, by Theorem \ref{thm3.2}, there exist
finitely many harmonic maps $\omega^{k}:S^{2}\rightarrow N,
k=1,2,...,K$, such that after selection of a subsequence, the
following holds:
\bee \label{4.14}\lim \limits_{n \rightarrow
\infty}E(u_{n},P_{n}^{\delta_{n}})= \sum_{k=1}^{K}E(\omega^{k})+\lim
\limits_{n \rightarrow \infty}|{\rm Re} \alpha_{n}|\cdot
\frac{\pi^{2}}{l_{n}}. \eee
Combining \eqref{4.7}, \eqref{4.13} and \eqref{4.14} gives
\bee \lim \limits_{n \rightarrow \infty}E(u_{n})&=&
E(u)+\sum_{i=1}^{I}\sum_{l=1}^{L_{i}}E(\sigma^{i,l}) +
\sum_{k=1}^{K}E(\omega^{k})+\lim \limits_{n \rightarrow \infty}|{\rm
Re} \alpha_{}|\cdot \frac{\pi^{2}}{l_{n}}. \nn \eee

Finally, we consider the general case $p>1$. By the thick-thin
decomposition of hyperbolic surfaces (c.f. \cite{Hu}, Lemma IV.4.1
and Proposition IV.4.2), both the short simple closed geodesics (of
lengths $< {\rm 2arcsinh(1)}$) and the corresponding ${\rm
arcsinh(1)}$-thin parts of the collars around them are pairwise
disjoint. Hence we can deal with the corresponding subcollars
separately, and the remaining proof is analogous to the simpler
case. This completes the proof. \hfill $\square$

\vskip 0.2cm

\noindent\emph{Proof of Theorem \ref{thm1.2}}. W.l.o.g., we assume
that $p=1$ and the limit $ \liminf \limits_{n \rightarrow
\infty}\sqrt{|{\rm Re} \alpha _{n}|}\cdot \frac{\pi^{2}}{l_{n}}$
exist in $[0,\infty]$. Then the results follow from applying Theorem
\ref{thm3.1} with domain cylinders $P_{n}^{\delta_{n}}$ as in the
proof of Theorem \ref{thm1.1}. \hfill $\square$

\vskip 0.2cm

\noindent\emph{Proof of Proposition \ref{pro1.1}}. By Theorem
\ref{thm1.1} and Proposition \ref{pro3.2}. \hfill $\square$

\vskip 0.2cm

\noindent{\bf Asymptotic behaviour.} For each $j$, the asymptotic
behaviour of the necks appearing near the $j$-th node is
characterized by $\{(\alpha_{n}^{j},l_{n}^{j})\}_{n=1}^{\infty}$,
namely
\bee \label{4.15} E^{j} \approx |{\rm Re} \alpha_{n}^{j}|\cdot
\frac{\pi^{2}}{l_{n}^{j}}, \qquad L^{j} \approx \sqrt{|{\rm Re}
\alpha _{n}^{j}|}\cdot \frac{\pi^{2}}{l_{n}^{j}}, \eee
where $E^{j}$ is the sum of the energies of the necks and $L^{j}$ is
the sum of the average lengths of the necks. In general, we have the
following four cases as $n \rightarrow \infty$:

\begin{itemize}
\item[{\rm (1)}]
$\quad E^{j} \rightarrow E_{0}, \hspace{14pt}  L^{j} \rightarrow
\infty,$

\item[{\rm (2)}]
$\quad E^{j} \rightarrow 0,   \qquad    L^{j} \rightarrow \infty,$

\item[{\rm (3)}]
$\quad E^{j} \rightarrow 0,  \qquad     L^{j} \rightarrow L_{0},$

\item[{\rm (4)}]
$\quad E^{j} \rightarrow 0,  \qquad     L^{j} \rightarrow 0.$
\end{itemize}
Here $E_{0}\in(0,\Lambda]$ and $L_{0}\in(0,\infty)$ are two
constants.

\begin{remark}\
\begin{itemize}
\item[(1)] If $u_{n}: (\Sigma,h_{n}) \rightarrow
(N,g)$ are conformal harmonic maps (i.e., minimal surfaces, in
particular, pseudo-holomorphic curves \cite{Gr}), i.e.,
$\Phi(u_{n})\equiv 0$, then it is easy to verify that $\alpha
_{n}^{j}\equiv0$, for all $n$ and $j$. It follows immediately that
$$\liminf \limits_{n \rightarrow \infty}|{\rm Re}
\alpha_{n}^{j}|\cdot \frac{\pi^{2}}{l_{n}^{j}}=0, \quad \liminf
\limits_{n \rightarrow \infty}\sqrt{|{\rm Re} \alpha _{n}^{j}|}\cdot
\frac{\pi^{2}}{l_{n}^{j}}=0, \quad \forall j=1,2,...,p.$$

\item[(2)] If, in addition, we assume that $u_{n}$ is an energy-minimizing sequence
in the same homotopy class, then, using a replacing argument from
\cite{CT}, one can show that in the limit the lengths of the necks
are all finite, i.e., $\liminf \limits_{n \rightarrow
\infty}\sqrt{|{\rm Re} \alpha _{n}^{j}|}\cdot
\frac{\pi^{2}}{l_{n}^{j}}< \infty $ for each $j$, which yields
$$\liminf \limits_{n \rightarrow \infty}|{\rm Re}
\alpha_{n}^{j}|\cdot \frac{\pi^{2}}{l_{n}^{j}}= 0, \quad \forall
j=1,2,...,p.$$


\item[(3)] When the domain surfaces are degenerating tori, the
problem is simpler for two reasons. Firstly, the moduli space of
complex structures on the torus is simple. Secondly, any holomorphic
quadratic differential on a torus is a constant. We refer to
\cite{Z} for more details. It is worth mentioning that Parker's
example \cite{Pa} illustrates the asymptotics that the necks become
longer and longer geodesics and carry a certain amount of energy.
Some modifications to his example can illustrate the four cases of
asymptotics mentioned, see \cite{Z}.

\end{itemize}
\end{remark}

\vskip 0.5cm

\noindent Miaomiao Zhu\\
Max Planck Institute for Mathematics in the Sciences\\
Inselstr.\ 22-26, D-04103 Leipzig, Germany\\
E-mail: Miaomiao.Zhu@mis.mpg.de

\end{document}